\documentclass[a4paper,11pt]{article}

\usepackage{braket}
\usepackage{listings}
\usepackage{amsmath}
\usepackage{amsfonts}
\usepackage{amsthm}
\usepackage{amssymb}
\usepackage{ascmac}
\usepackage{fancybox}

\usepackage{newtxtext}

\usepackage{titlesec}
\usepackage{here}
\usepackage{bm}
\usepackage{caption}
\usepackage{enumerate}
\usepackage{pdfpages}

\usepackage{siunitx}
\usepackage{url}

\usepackage{a4wide}

\usepackage{graphicx}

\usepackage[linktocpage, colorlinks=true]{hyperref}
\usepackage[capitalise]{cleveref}
\usepackage{caption}

\theoremstyle{definition}
\newtheorem{theorem}{Theorem}[section]

\newtheorem{lemma}[theorem]{Lemma}

\makeatletter
\@addtoreset{figure}{section}
\@addtoreset{table}{section}
\@addtoreset{equation}{section}
\makeatother

\usepackage{docmute}

\begin{document}
\title{Maximum Cluster Diameter in Non-Critical Bond Percolation}
\author{Kaito Kobayashi}
\date{}

\maketitle

\begin{abstract}
  In this paper, we study independent (Bernoulli) bond percolation in dimensions $d \ge 2$, focusing on the maximum diameter of finite clusters in the non-critical regime ($p\neq p_c$). We prove that the maximum diameter $R_n$ satisfies $R_n / \log n \to \varkappa(p)$ almost surely, where $\varkappa(p)$ is determined by the exponential decay rate $\xi(p)$ of $P_p(0 \leftrightarrow \partial B_n, |\mathcal C_0|<\infty)$. Furthermore, we establish a large deviation principle for the event $\{R_n > \rho\log n\}$ for $\rho > \varkappa (p)$. Finally, we consider the asymptotics of the number of vertices in clusters with large diameters.
\end{abstract}


\section{Introduction}
Percolation is a fundamental stochastic model introduced by Broadbent and Hammersley in 1957 to formalize the permeation of liquids and gases through porous media \cite{BH1957}. In recent years, its applications have expanded across various fields, ranging from phase transitions in mathematical physics to the propagation of viruses; see, e.g., \cite{SCA1982, SS, GGP2}. A central feature of the model involves random subgraphs of a given domain, where the connectivity properties change drastically as a density parameter varies. This phenomenon is visually illustrated in Figure \ref{fig:phase_transition}.

Fundamental questions in percolation theory concern the existence and uniqueness of the infinite cluster, the critical value, and the decay of correlation functions. While these global properties have been extensively studied, the geometric characteristics of finite connected components, particularly their shape and spatial extent, are equally important for understanding the underlying graph structure.

Historically, much interest has focused on the size of finite clusters. The tail behavior of the radius of an open cluster, in the sense of Grimmett (i.e., the event $0\leftrightarrow \partial B_n$), was first estimated by Hammersley \cite{JMH1957}, and since then, numerous studies have refined these estimates; see, e.g., \cite{KS1978, ADS1980}. More complex geometric properties have been thoroughly investigated at the critical value $p_c (d)$, where scaling exponents such as the fractal dimension and the chemical distance exponent play a central role; see, e.g., \cite{K1987, HvH2007, HK2010, KB2016}.

However, compared to the critical regime, the geometry of clusters in the non-critical regime ($p \neq p_c$) is much less understood. 
Although the expected size and the tail behavior of finite clusters are well understood, sharp asymptotics for geometric quantities, most notably the maximum diameter of finite clusters measured in the $\ell^\infty$ metric, remain largely unknown in both the subcritical and supercritical regimes.
This paper addresses this gap.

\subsection{The model}
We consider $d \ge 2$ bond percolation on the $d$-dimensional integer lattice $(\mathbb{Z}^d, \mathbb{E}^d)$, where $\mathbb{Z}^d$ is the set of vertices and $ \mathbb{E}^d$ is the set of edges.
Each edge is declared open with probability $p$ and closed with probability $1-p$, independently of the other edges.
The resulting configuration of the bonds is denoted by $\omega \in \Omega$, where $\Omega = \{0, 1\}^{\mathbb{E}^d}$ is the configuration space. 
A bond $b$ is called open if $\omega(b) = 1$, and closed if $\omega(b) = 0$.
The corresponding probability measure on $\{ 0,1 \}^{\mathbb{E}^d}$ is denoted by $P_p$, and the corresponding expectation operator is denoted by $E_p$. 
For $x \in \mathbb{Z}^d$, we denote by $\mathcal{C}_x = \mathcal{C}_x (\omega)$ the cluster of $x$. 
We write $x \leftrightarrow y$ to denote the event that the vertices $x$ and $y$ belong to the same cluster.
Throughout this paper, we assume that $p \neq p_c$; that is, $p$ is strictly different from the critical probability
\[
p_c (d) := \sup \{ p : \theta (p) = 0\},
\]
where $\theta (p)$ is the probability that the origin belongs to an infinite cluster.
\begin{figure} 
  \centering
  \begin{minipage}{0.45\columnwidth}
  \captionsetup{labelformat=empty, skip=0pt}
    \centering
    \includegraphics[width=0.9\columnwidth]{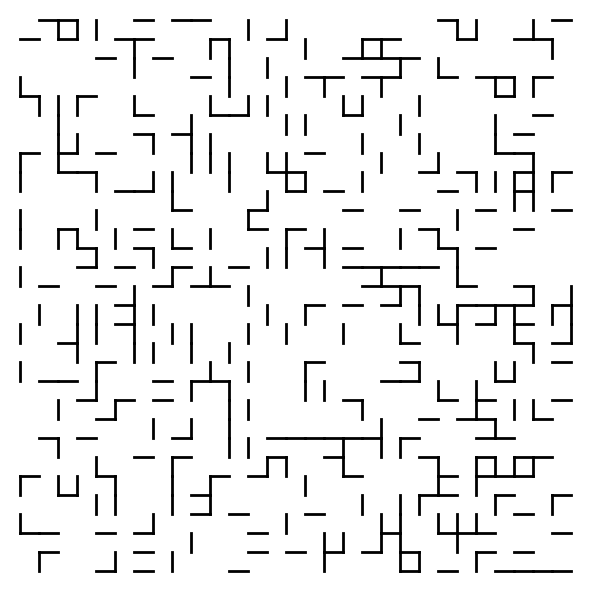}
    \caption{(a) \, $p=0.25$}
  \end{minipage}
  \hfill
  \begin{minipage}{0.45\columnwidth}
  \captionsetup{labelformat=empty, skip=0pt}
    \centering
    \includegraphics[width=0.9\columnwidth]{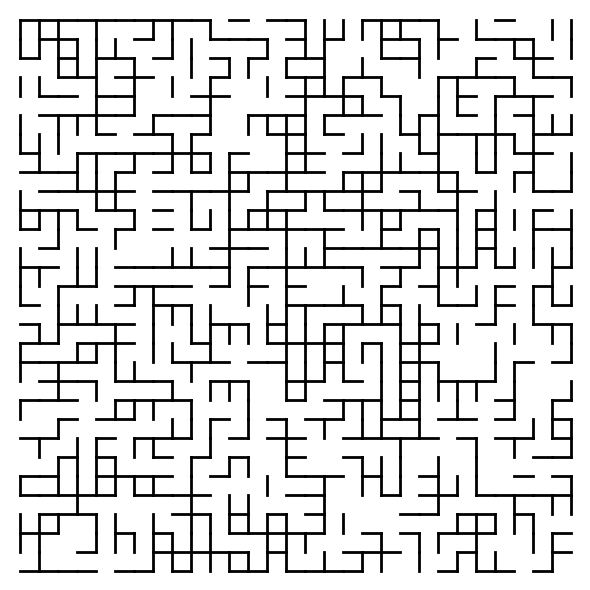}
    \caption{(b) \, $p=0.49$}
  \end{minipage}
  \vspace{5pt}
  \begin{minipage}{0.45\columnwidth}
  \captionsetup{labelformat=empty, skip=0pt}
    \centering
    \includegraphics[width=0.9\columnwidth]{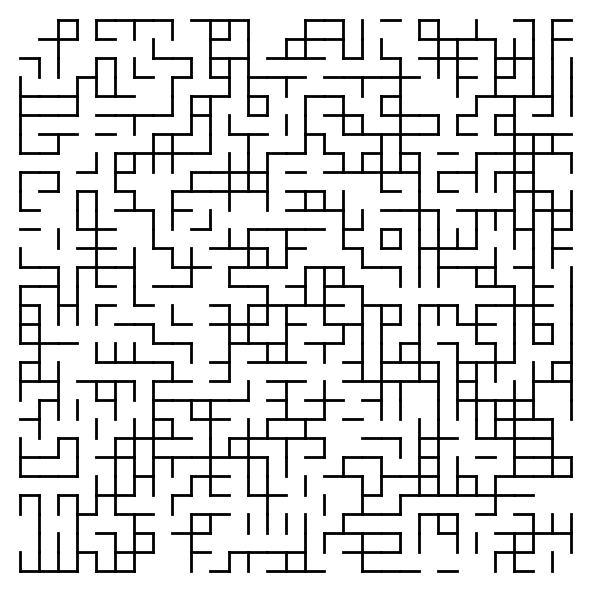}
    \caption{(c) \, $p=0.51$}
  \end{minipage}
  \hfill
  \begin{minipage}{0.45\columnwidth}
  \captionsetup{labelformat=empty, skip=0pt}
    \centering
    \includegraphics[width=0.9\columnwidth]{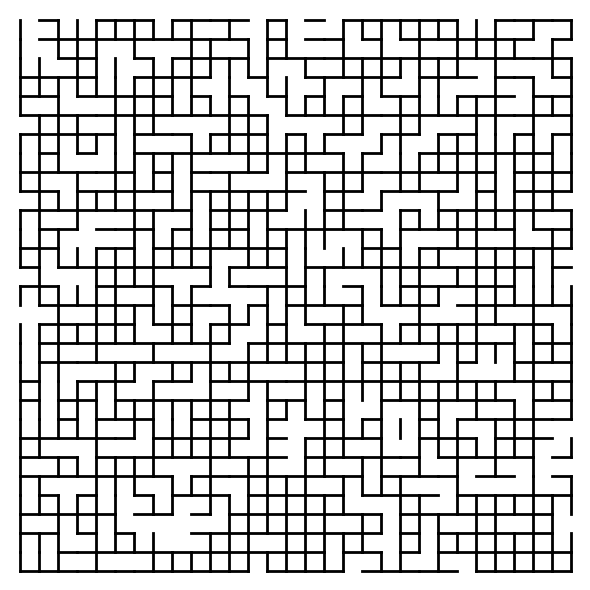}
    \caption{(d) \, $p=0.75$}
  \end{minipage}
  \caption{Realizations of bond percolation on a $30 \times 30$ section of the square lattice for four different values of $p$.}
  \label{fig:phase_transition}
\end{figure}
\newpage
Broadbent and Hammersley \cite{BH1957} established the existence of a critical probability $p_c$, setting the stage for the proof that $p_c (d) \in (0,1)$ for all $d \ge 2$, which was later rigorously completed.
This value marks the boundary where the model switches from having no infinite cluster to having exactly one.
In the two-dimensional case, Kesten \cite{KESTEN1980} proved that $p_c (2) = 1/2$.
As illustrated in Figure \ref{fig:phase_transition} (b) and (c), open paths connect the left and right sides of the box when $p=0.51$, whereas no such path exists when $p=0.49$.
Furthermore, Aizenman, Kesten, and Newman \cite{AKN1987} proved that the infinite cluster is unique whenever it exists.
Later, an alternative, significantly simpler and more elegant proof of this uniqueness result was provided by Burton and Keane \cite{BK1989}.
Although the specific values of $p_c (d)$ for $d \ge 3$ remain an important unsolved problem, the field has seen remarkable progress \cite{AGKSZ2014, DC2018, G2023, M2017, BS1996}.

A key feature of the non-critical regime ($p \neq p_c$) is the rapid spatial decay of connectivities within finite clusters.
In the subcritical phase ($p < p_c$), the radius of a finite open cluster decays exponentially, a result rigorously established by Menshikov \cite{M1986} and Aizenman and Barsky \cite{AB1987}.
In the supercritical phase ($p > p_c$), while a unique infinite cluster exists, the distribution of the size of the remaining finite clusters exhibits similar exponential decay properties \cite{CCGKS1989}.

The main objective of this paper is to investigate the asymptotic behavior of the maximum diameter of finite clusters.
To make this precise, we introduce the necessary geometric definitions and notation in the following subsection.

\begin{figure}
    \centering
    \includegraphics[width=0.4\linewidth]{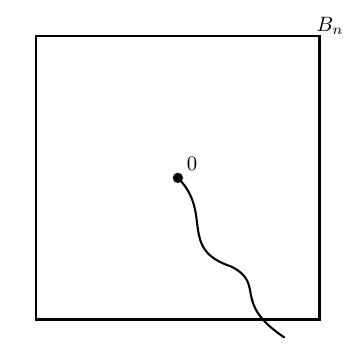}
    \caption{A sketch of the event $\{0 \leftrightarrow \partial B_n , |\mathcal{C}_0| < \infty\}$: the origin belongs to a finite cluster which intersects $\partial B_n$.}
    \label{fig:event:0toBn}
\end{figure}

\subsection{The maximum diameter of finite open clusters}
In this section, we introduce the definitions and notation concerning the diameter of clusters.

For $x \in \mathbb{Z}^d$, we write $x_i$ for the $i$th coordinate of $x$. We denote by $B_n$ the box with side length $2n + 1$ and center at the origin: 
\[
B_n := \{ x \in \mathbb{Z}^d : \|x\|_{\infty} \le n \},
\]
where $\|x\|_{\infty} = \max \{|x_i| : 1 \le i \le d \}$.
We write $\partial B_n$ for the boundary of $B_n$: 
\[
\partial B_n := \{ x \in \mathbb{Z}^d : \|x\|_{\infty} = n \}.
\]
As usual, $P_p (0 \leftrightarrow \partial B_n)$ is the probability that there exists an open path joining the origin to some vertex in $\partial B_n$; see Figure \ref{fig:event:0toBn}.
We refer to $P_p(0 \leftrightarrow \partial B_n)$ as the one-arm probability.

 We define the maximum diameter of finite clusters by
\begin{align}\label{def:R_n^fb}
R_n
&:= \max\Bigl(\{ \mathrm{diam}(\mathcal{C}_x(\omega)) :
x\in B_n, |\mathcal{C}_x(\omega)|<\infty \}\cup\{0\}\Bigr) \notag\\
&= \max\Bigl(\{ \max_{1\le i\le d} D_i(\mathcal{C}_x(\omega)) :
x\in B_n, |\mathcal{C}_x(\omega)|<\infty \}\cup\{0\}\Bigr).
\end{align}
where for $i = 1,2,\ldots,d$,
\[
D_i (\mathcal{C}) := \max_{y \in \mathcal{C}} y_i - \min_{y \in \mathcal{C}} y_i.
\]

Specifically, $D_i$ is defined as the maximum width of any cluster $\mathcal{C}$ in the $i$th coordinate direction, and the maximum diameter of a finite cluster is defined as the maximum value of $D_i$.

We note that the distribution of the diameter plays an essential role throughout the proofs.
To evaluate the asymptotic behavior of $R_n$, it is necessary to consider how the one-arm probability decays with respect to distance. 
We now state the necessary results concerning the distribution of the diameter. 
Since the distribution law fundamentally differs between the subcritical ($p < p_c$) and supercritical ($p > p_c$) regimes, we treat these two cases separately.

For $0 < p < p_c$, it is shown in \cite[Theorem 6.10]{GGP2} that
\begin{align} \label{def:varphi,Grimmett:6.10}
    \varphi (p) = \lim_{n \to \infty} \left\{ - \frac{1}{n} \log P_p ( 0 \leftrightarrow \partial B_n )\right\} ,
\end{align}
exists, and satisfies $\varphi (p) > 0$ if $p \in (0,p_c)$ \cite[Theorem 6.14]{GGP2}. 
Furthermore, by \cite[(6.11)]{GGP2}, there exists a constant $a>0$, independent of $p$, such that, for all $n$,
\begin{align} \label{Grimmett:6.11}
    P_p (0 \leftrightarrow \partial B_n) \le an^{d-1} e^{-n \varphi (p)}.
\end{align}
We now turn to supercritical results.
When $p > p_c$, the probability of the origin belonging to the infinite cluster is positive $\bigl(P_p(|\mathcal{C}_0|=\infty) > 0 \bigr)$. 
Therefore, instead of evaluating the exponential decay of the event $\{0 \leftrightarrow \partial B_n\}$, we must specifically consider the probability under the condition that the origin belongs to a finite cluster.
For $p > p_c$, it is shown in \cite[Theorem 8.18]{GGP2} that
\begin{align} \label{def:sigma,Grimmett:8.18}
    \sigma(p) = \lim_{n \to \infty} \left\{ -\frac{1}{n} \log P_p (0 \leftrightarrow \partial B_n, |\mathcal{C}_0| < \infty) \right\}
\end{align}
exists and satisfies $\sigma (p) < \infty$.
Furthermore, there exists a constant $A(p,d)$ which is finite for $d \ge 2, 0 < p < 1$, such that, for all $n$,
\begin{align} \label{Grimmett:8.20}
    P_p (0 \leftrightarrow \partial B_n, |\mathcal{C}_0| < \infty) \leq A(p,d) n^d e^{-n\sigma (p)}.
\end{align}
In fact, by \cite[Theorem 8.21]{GGP2}, we have that $\sigma (p) > 0$ if $p > p_c$.

For $p \neq p_c$, we define $\xi (p)$ as follows:
\[
\xi (p) := \lim_{n \to \infty} \left\{ -\frac{1}{n} \log P_p (0 \leftrightarrow \partial B_n, |\mathcal{C}_0| < \infty) \right\}.
\]
In particular, in view of \eqref{def:varphi,Grimmett:6.10} and \eqref{def:sigma,Grimmett:8.18}, we have 
\begin{align*}
    \xi(p)=
    \begin{cases}
        \varphi (p), &0 < p <p_c,\\
        \sigma (p), &p>p_c.
    \end{cases}
\end{align*}
Note that the decay rate of $P_p (0 \leftrightarrow \partial B_n, |\mathcal{C}_0| < \infty)$ matches the decay rate of $P_p (\mathrm{diam}(\mathcal{C}_0) = n)$. 
Moreover, the identity below was proved in \cite[(8.37)]{GGP2}:
\begin{align} \label{def:xi}
   \xi (p) =  \lim_{n \to \infty} \left\{ - \frac{1}{n} \log P_p \bigl(\textnormal{diam}(\mathcal{C}_0) = n\bigr) \right\}>0.
\end{align}
We note that 
\[
P_p \bigl(\textnormal{diam}(\mathcal{C}_0) \ge n , |\mathcal{C}_0| < \infty \bigr) = \sum_{m=n}^{\infty} P_p \bigl(\textnormal{diam}(\mathcal{C}_0) =m \bigr),
\]
so that the limit $\xi (p)$ in \eqref{def:xi} satisfies
\begin{align} \label{def:xi'}
    \xi (p) = \lim_{n \to \infty} \left\{ - \frac{1}{n} \log P_p \bigl(\textnormal{diam}(\mathcal{C}_0) \ge n , |\mathcal{C}_0| < \infty \bigr) \right\}.
\end{align}
Thus, by \eqref{Grimmett:6.11} and \eqref{Grimmett:8.20}, there exists a constant $L$ for $p \neq p_c$, such that
\begin{align} \label{ineq:Pp(diam)<=Ke^-xi}
    P_p \bigl(\textnormal{diam}(\mathcal{C}_0) \ge n , |\mathcal{C}_0| < \infty \bigr) \le L n^d e^{-n\xi (p)}.
\end{align}
Furthermore, we define the constant $\varkappa(p)$ as follows:
\[
\varkappa(p) := \frac{d}{\xi(p)}.
\]
We are now able to state our results.

\section{Main Results}
The next theorem shows the typical size of $R_n$ in more detail and establishes almost sure convergence.
\begin{theorem} \label{thm:LLN:Rn}
For $p \neq p_c$,
\begin{align*}
    \frac{R_n}{\log n} \to \varkappa(p) \qquad \text{a.s.}
\end{align*}
\end{theorem}
We now define and discuss the boundary conditions for $R_n$.
In \eqref{def:R_n^fb}, we have taken the maximum diameter under the free boundary condition, so that we can write $R_n = R_n^{(\mathrm{fb})}$. 
Alternatively, we consider $R_n$ under the zero boundary condition, i.e., 
$R_n^{(\mathrm{zb})}$ is defined by restricting the configuration $\omega$ to the box $B_n$.

\begin{figure} 
  \centering
  \begin{minipage}[b]{0.48\columnwidth}
    \captionsetup{labelformat=empty, skip=0pt}
    \centering
    \includegraphics[width=1\columnwidth]{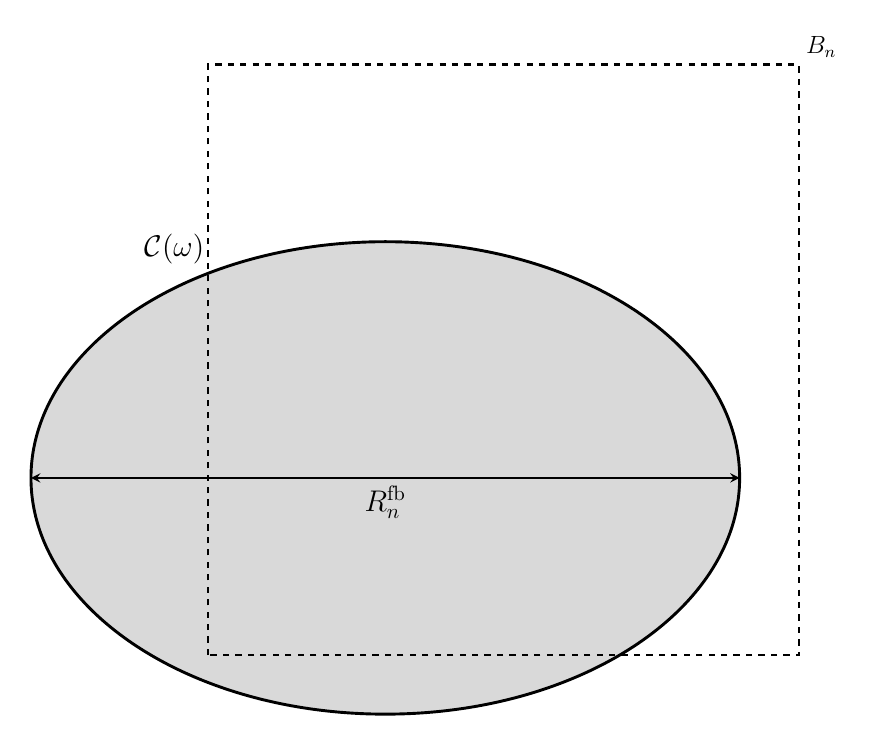}
  \end{minipage}
      \raisebox{5mm}{
  \begin{minipage}[b]{0.48\columnwidth}
    \captionsetup{labelformat=empty, skip=0pt}
    \centering
    \includegraphics[width=0.8\columnwidth]{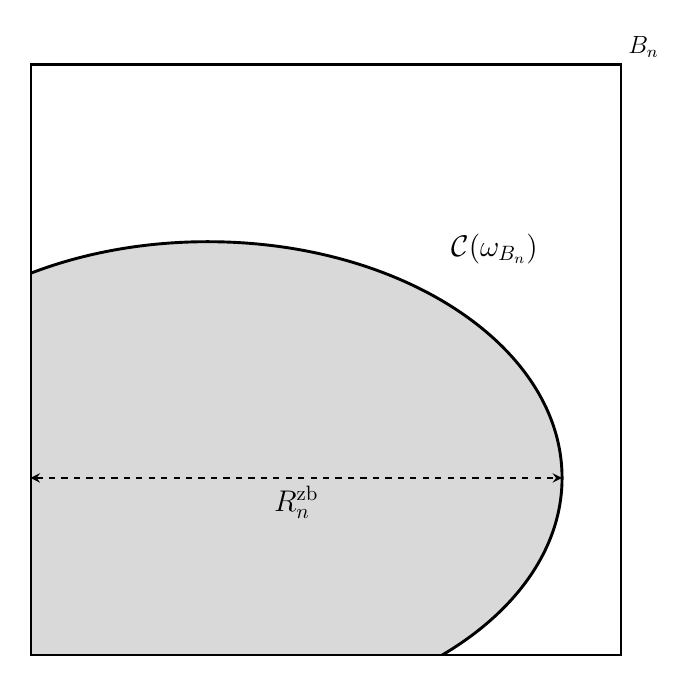}
  \end{minipage}
  }
  \caption{Visualizing the difference between $R_n^{(\mathrm{fb})}$ and $R_n^{(\mathrm{zb})}$. The probability that they differ tends to $0$ as $n\to\infty$ (Theorem \ref{thm:Rnfb=Rnzb}).}
  \label{fig:differenceRfbRzb}
\end{figure}

We define the restriction $\omega_{B_n}$ of $\omega$ to $B_n$ as follows:
\begin{align*}
\omega_{B_{n}}(b)= 
\begin{cases}
\omega(b) & \text{if } b\subset  B_{n},  \\[6pt]
0 & \text{otherwise}.
\end{cases}
\end{align*}
That is, $\omega_{B_n}$ is the configuration where all bonds existing outside $B_n$ are forced to be closed. 
By introducing this restriction, it becomes possible to capture local geometric features, such as cluster diameter and size, within the finite region $B_n$.
We then define $R_n^{(\mathrm{zb})}$, the maximum cluster diameter under this zero boundary condition, as
\[
R_n^{(\mathrm{zb})}
:=\max\Bigl(\{\mathrm{diam}(\mathcal{C}_x(\omega_{B_n})): x\in B_n, |\mathcal{C}_x(\omega)|<\infty\}\cup\{0\}\Bigr).
\]
We restrict to vertices that belong to finite clusters in $\omega$ in order to exclude contributions from the infinite cluster in the supercritical regime.
Figure \ref{fig:differenceRfbRzb} illustrates the distinction between $R_n^{(\mathrm{fb})}$ and $R_n^{(\mathrm{zb})}$.
Despite this difference, the following theorem shows that they are asymptotically equivalent.

\begin{theorem} \label{thm:Rnfb=Rnzb}
For $p \neq p_c$,
\begin{align*}
 \lim_{n\to \infty}   P_p \bigl(R_n^{(\mathrm{zb})} \neq R_n^{(\mathrm{fb})} \bigr) = 0.
\end{align*}
\end{theorem}

Theorem \ref{thm:LLN:Rn} shows that $R^{({\mathrm{fb}} )}_n / \log n$ converges to $\varkappa(p)$ almost surely as $n \to \infty$. The next theorem gives a large deviation principle describing the probability that $R_n$ deviates from its typical value.

\begin{theorem} \label{thm:LDP:Rn}
For any $p \neq p_c$ and $\rho > \varkappa(p)$, the limit 
\[
\varsigma(\rho) =
\lim_{n \to \infty} -\frac{1}{\log n} \log P_p\big(R_n > \rho \log n \big),
\]
exists and satisfies $\varsigma (\rho) > 0$.
Furthermore, $\varsigma(\rho)$ is given by
\[
\varsigma(\rho) = \xi(p) \rho - d.
\]
\end{theorem}

Finally, our result connects the geometric characteristics established by Theorem \ref{thm:LDP:Rn} to the spatial extent of these rare events. The following theorem describes the asymptotic behavior of the number of vertices belonging to clusters whose diameter exceeds $\rho \log n \; (\rho < \varkappa(p))$, using the exponent $d - \xi(p)\rho$ appearing in the large deviation principle.
\begin{theorem} \label{thm:ProbCon:Sn}
Let $d \ge 2$ and $p \neq p_c$. For any $0 < \rho < \varkappa(p)$, let $S_n(\rho)$ be the number of vertices in $B_n$ that belong to a finite cluster whose diameter exceeds $\rho \log n$, defined as
\[
S_n(\rho) := |\{x\in B_n: ~\rho \log n<\mathrm{diam}(\mathcal{C}_x)<\infty\}|.
\]
Then, the following holds:
\[
\frac{S_n(\rho)}{E_p[S_n(\rho)]} \xrightarrow{P_p} 1 \quad \text{as } n \to \infty.
\]
Furthermore, for any $\varepsilon>0$, when $n$ is sufficiently large,
\[
n^{d - \xi(p)\rho -\varepsilon} \leq E_p[S_n(\rho)]\leq n^{d - \xi(p)\rho +\varepsilon}.
\]
\end{theorem}
\subsection{Related work}
The asymptotic properties of the maximal volume of finite clusters in the non-critical regime were investigated by van der Hofstad and Redig \cite{HR2006}, who studied the typical behavior of the maximum cluster volume within a finite box and established its logarithmic scaling.
In this work, we extend their framework to analyze the diameter of clusters, and this extension is a key ingredient in the proofs of Theorems \ref{thm:LLN:Rn} and \ref{thm:Rnfb=Rnzb}.
Although the definition of diameter differs from that of the volume, the underlying probabilistic structure regarding the exponential decay allows us to adapt the techniques of \cite{HR2006}. In addition, we establish a large deviation principle for the diameter and the asymptotics of the number of vertices in clusters with large diameters.

We also remark that our results may be extended to diameters measured in the chemical (graph) distance, rather than the $\ell^\infty$-norm used here. This extension is not immediate, but it is supported by the recent work of Dembin and Nakajima \cite{DN2022}, which established the existence of a specific rate function for the chemical distance. Their result indicates that, although additional work is required, one can expect that the maximal chemical diameter among finite clusters is characterized by this rate function. 
\section{Proofs}
In this section, we prove the main results stated in Section 2.

\subsection{Preliminaries}
This section provides the necessary lemmas for the proofs.
\begin{lemma}[Borel–Cantelli] \label{lem:BorelCantelli}
Let $\{E_n\}_{n=1}^{\infty}$ be a sequence of events. If $\sum_{n=1}^{\infty} P(E_n) < \infty$,
then
\[
P \left(\limsup_{n \to \infty} E_n \right) = 0.
\]
\end{lemma}

\begin{lemma}[Boole's inequality] \label{lem:unionbound}
For any sequence of events $\{A_i\}$, we have
\[ P\left(\bigcup_i A_i\right) \le \sum_i P(A_i). \]
\end{lemma}

\subsection{Proof of Theorem \ref{thm:LLN:Rn} and \ref{thm:Rnfb=Rnzb}}
We adapt the methods used by van der Hofstad and Redig \cite[Theorems 1.3-1.5]{HR2006}.
\begin{proof}[\textnormal{\textbf{Proof of Theorem \ref{thm:LLN:Rn}}}]
The main ingredient is the following lemma:
\begin{lemma} \label{lem:forpfofthm}
For any $\varepsilon > 0$, there exists $\kappa > 0$ such that as $n \to \infty$,
\begin{align*}
    P_p \left(\left|\frac{R_n}{\log n} -\varkappa(p)\right|>\varepsilon\right) < n^{-\kappa}.
\end{align*}
\end{lemma}

Before proving Lemma \ref{lem:forpfofthm}, we first complete the proof of Theorem \ref{thm:LLN:Rn} assuming Lemma \ref{lem:forpfofthm}.

Let $n_k = 2^k$.
As a consequence of Lemma \ref{lem:forpfofthm}, and the fact that for every
$\kappa > 0$,
\[
n^{-\kappa}_k = 2^{-\kappa k}
\]
is summable in $k$, so that
\begin{align*}
    P_p \left( \limsup_{k \to \infty} \left|\frac{R_{n_k}}{\log n_k} -\varkappa(p)\right|>\varepsilon\ \right) = 0
\end{align*}
by Lemma \ref{lem:BorelCantelli}.
Thus,
\begin{align*}
    \frac{R_{n_k}}{\log n_k} \to \varkappa(p)  \qquad \text{a.s.}
\end{align*}
Moreover, $n \mapsto R_n$ is non-decreasing with probability $1$.
Therefore, for any $n_k < n \leq n_{k+1}$ we can bound
\begin{align} \label{boundineq:RnkRnRnk+1}
    \frac{R_{n_k}}{\log (n_k)} \frac{\log (n_k)}{\log (n_{k+1})} \le \frac{R_n}{\log n} \le \frac{R_{n_{k+1}}}{\log (n_{k+1})} \frac{\log (n_{k+1})}{\log (n_{k})}.
\end{align}
As $n \to \infty$, also $n_k, n_{k+1} \to \infty$.
Thus, $\frac{R_{n_k}}{\log (n_k)}$ and $\frac{R_{n_{k+1}}}{\log (n_{k+1})}$ converge to $\varkappa(p)$ with probability $1$.
Furthermore, 
\begin{align*}
\lim_{k \to \infty} \frac{\log (n_{k+1})}{\log (n_{k})}=
\lim_{k \to \infty} \frac{\log (n_{k})}{\log (n_{k+1})} = \lim_{k \to \infty} \frac{k+1}{k} = 1,
\end{align*}
so that both upper and lower bounds in \eqref{boundineq:RnkRnRnk+1} converge to $\varkappa(p)$ almost surely.
This completes the proof of Theorem \ref{thm:LLN:Rn}.
\end{proof}

\begin{proof}[\textnormal{\textbf{Proof of Lemma \ref{lem:forpfofthm}}}]
Let $\varepsilon > 0$. Without loss of generality, we assume that $\varepsilon\in (0,\varkappa(p))$. We set $ \delta = \frac{\xi(p)\,\varepsilon}{2(\varkappa(p)+\varepsilon)}$ and $\kappa' = \xi (p) \varepsilon / 4 > 0$, using \eqref{def:xi'} and $\xi (p)=d/\varkappa(p)$.
Then, the following holds:
\begin{align} \label{Ppineq:diam>,<ndk'}
    P_p \left( \frac{\mathrm{diam}(\mathcal{C}_0)}{\log n} > \varkappa(p) + \varepsilon \,;\,\mathrm{diam}(\mathcal{C}_0)<\infty\right) \notag
    &\le \exp \left\{ -(\xi (p) -\delta)(\varkappa(p) + \varepsilon)\log n \right\}\\ \notag
    &= n^{-(d+\xi (p) \varepsilon -\delta ( \varkappa(p) + \varepsilon))}\\ 
    &\leq n^{-d - \kappa' }.
\end{align}
Moreover,  setting $ \delta' = \frac{\xi(p)\varepsilon}{2(\varkappa(p)-\varepsilon)}$ and using the same $\kappa' = \xi (p) \varepsilon/4>0$, for sufficiently large $n$, the following holds:
\begin{align}\label{Ppineq:diam<,<ndkc}
&P_p\!\left(
  \frac{\mathrm{diam}(\mathcal{C}_0)}{\log n}
  \notin (\varkappa(p)-\varepsilon,2\varkappa(p)] \,;\,
  \mathrm{diam}(\mathcal{C}_0)<\infty
\right) \notag\\
&\le P_p(|\mathcal{C}_0|<\infty)
   - P_p\!\left(
       \frac{\mathrm{diam}(\mathcal{C}_0)}{\log n}>\varkappa(p)-\varepsilon \,;\,
       \mathrm{diam}(\mathcal{C}_0)<\infty
     \right) \notag\\
&\hphantom{= P_p(|\mathcal{C}_0|<\infty)\;}
   + P_p\!\left(
       \frac{\mathrm{diam}(\mathcal{C}_0)}{\log n}>2\varkappa(p) \,;\,
       \mathrm{diam}(\mathcal{C}_0)<\infty
     \right) \notag\\
&\le 1 - \exp\!\left\{
        -\left(\xi(p)+\delta'\right)\left(\varkappa(p)-\varepsilon\right)\log n
      \right\}
      + n^{-d-\kappa'} \notag\\
&\le 1 - \frac12\,n^{-d+2\kappa'}.
\end{align}

Fix $\varepsilon \in (0,\varkappa(p))$. We will prove that there exists $\kappa > 0$ such that
\begin{align}\label{Ppineq:Rn>,<nk}
    P_p \left( \frac{R_n}{\log n} > \varkappa(p) + \varepsilon\right) \leq n^{ - \kappa}, 
\end{align}
and
\begin{align}\label{Ppineq:Rn<,<nk}
    P_p \left( \frac{R_n}{\log n} < \varkappa(p) - \varepsilon\right) \leq n^{ - \kappa}. 
\end{align}

To prove \eqref{Ppineq:Rn>,<nk}, by Lemma \ref{lem:unionbound} and \eqref{Ppineq:diam>,<ndk'},
\begin{align*}
    P_p \left( \frac{R_n}{\log n} > \varkappa(p) + \varepsilon\right) \notag
    &= P_p \left( \bigcup_{x \in B_n} \bigl\{ (\varkappa(p) + \varepsilon)(\log n) < \mathrm{diam}(\mathcal{C}_x) < \infty  \bigr\}\right) \notag \\
    & \leq \sum_{x \in B_n} P_p \bigl((\varkappa(p) + \varepsilon)(\log n) <\mathrm{diam}(\mathcal{C}_x) < \infty \bigr) \notag \\
    & \leq |B_n|n^{- d - \kappa'}.
\end{align*}
We can estimate the right-hand side as follows for sufficiently large $n$: 
\begin{align*}
        & \left( {2n + 1} \right)^d n^{- d - \kappa'} 
        \leq 3^d n^{- \kappa'} \leq n^{- \kappa'/2}.
\end{align*}
Thus, setting $\kappa=\kappa'/2$, we obtain \eqref{Ppineq:Rn>,<nk}.

To prove \eqref{Ppineq:Rn<,<nk}, let $\varepsilon\in(0,\varkappa(p))$, 
\[
A_n = (K_n \mathbb{Z})^d \cap B_n,
\]
and
\[
K_n = \lceil 8 \varkappa(p)\log n \rceil.
\] We use the fact that the events 
\[
\bigl\{ \mathrm{diam}(\mathcal{C}_x) \notin (
    (\varkappa(p) - \varepsilon)\log n,2\varkappa(p) \log n] \bigr\}_{x \in A_n}
\]
are independent, since the event $\{\mathrm{diam}(\mathcal{C}_x) \in (
    (\varkappa(p) - \varepsilon)\log n,2\varkappa(p) \log n]\}$ depends only on the states of edges inside $x+B_{\lceil 2\varkappa(p)\log n\rceil+1}$.
By independence,
\begin{align*}
    P_p \left( \frac{R_n}{\log n} \leq \varkappa(p) - \varepsilon\right) \notag
    &\leq P_p \left( \bigcap_{x \in A_n} \bigl\{ \mathrm{diam}(\mathcal{C}_x) \notin (
    (\varkappa(p) - \varepsilon)\log n,2 \varkappa(p) \log n] \bigr\}\right) \notag \\
    &=
    \prod_{x \in A_n} P_p \bigl(\mathrm{diam}({\mathcal{C}_x}) \notin ( (\varkappa(p) - \varepsilon)\log n,2\varkappa(p) \log n] \bigr) \notag \\
    &\leq P_p \bigl(\mathrm{diam}(\mathcal{C}_0) \notin ( (\varkappa(p) - \varepsilon)\log n, 2\varkappa(p) \log n] \bigr)^{|A_n|}.
\end{align*}
Next, by \eqref{Ppineq:diam<,<ndkc} and the fact that
\begin{align*}
    |A_n| \geq \left( \frac{n}{ 8 \varkappa(p) \log n}\right)^d,
\end{align*}
with the same $\kappa'=\xi(p)\varepsilon/4>0$, 
\begin{align*}
    P_p \left( \frac{R_n}{\log n} \leq \varkappa(p) - \varepsilon\right) 
    &\leq (1 - \frac{1}{2}n^{-d + 2\kappa'})^{|A_n|}\\
    &\leq \left(1 - \frac{1}{2}n^{-d+2\kappa'}\right)^{\left(\frac{n}{8 \varkappa(p) \log n }\right)^d} \\
    &\leq \left\{\exp\left(- \frac{1}{2}n^{-d+2\kappa'}\right)\right\} ^{\left(\frac{n}{8 \varkappa(p) \log n} \right)^d} \\
    &= \exp\left\{-\frac{ \frac{1}{2}n^{2\kappa'}}{(8 \varkappa(p) \log n)^{d}}\right\}.
\end{align*}
Therefore, we obtain \eqref{Ppineq:Rn<,<nk}.
The proof is complete.
\end{proof}

\begin{proof}[\textnormal{\textbf{Proof of Theorem \ref{thm:Rnfb=Rnzb}}}]
We use the fact that the event $\{R_n^{(\mathrm{zb})} \neq R_n^{(\mathrm{fb})}\}$ is contained in the event that there exists some vertex $x \in \partial B_n$ satisfying $\mathrm{diam}(\mathcal{C}_x) \ge R_n^{\mathrm{(fb)}}$. By Theorem \ref{thm:LLN:Rn}, for any $\varepsilon > 0$, $R_n^{(\mathrm{fb})} > (\varkappa(p) - \varepsilon)\log n$ holds with high probability. Thus, by Lemma \ref{lem:unionbound}, the probability can be bounded above as follows: 
\begin{align*}
    P_p(R_n^{(\mathrm{zb})} \neq R_n^{(\mathrm{fb})})
    &\le P_p \left( \bigcup_{x \in \partial B_n} \left\{R_n^{(\mathrm{fb})} \le\mathrm{diam}(\mathcal{C}_x) < \infty \right\}\right) \\
    &\le P_p \left( \bigcup_{x \in \partial B_n} \left\{(\varkappa(p) - \varepsilon)\log n \le \mathrm{diam}(\mathcal{C}_x) < \infty\right\} \right) +o(1)\\
    &\le \sum_{x \in \partial B_n} P_p \bigl( (\varkappa(p) - \varepsilon)\log n \le \mathrm{diam}(\mathcal{C}_x) < \infty \bigr) +o(1).
\end{align*}
Furthermore, by setting $\varepsilon = 1/(2\xi (p))$ and combining this with \eqref{def:xi'}, the right-hand side is bounded above by:
\begin{align*}
    &2d(2n+1)^{d-1} P_p\bigl((\varkappa(p) - \varepsilon)\log n \le \mathrm{diam}(\mathcal{C}_x) < \infty \bigr) \\
    &\le 2d 3^{d-1} n^{d-1} n^{-d +\xi (p)\varepsilon + o(1)} \\
    &= (2d 3^{d-1}) n^{-1+\xi (p)\varepsilon +o(1)} \to 0 \quad\text{ as } n \to \infty,
\end{align*}
where we have used $| \partial B_n| \leq 2d(2n+1)^{d-1}$.
\end{proof}

\subsection{Proof of Theorem \ref{thm:LDP:Rn} and \ref{thm:ProbCon:Sn}}
\begin{proof}[\textnormal{\textbf{Proof of Theorem \ref{thm:LDP:Rn}}}]
To prove the existence of the rate function, we proceed by proving the upper and lower bounds separately.
We first establish the upper bound. 
Using Lemma \ref{lem:unionbound} and \eqref{ineq:Pp(diam)<=Ke^-xi} with some constant $L$, the probability is bounded  as follows: 
\begin{align*}
P_p(R_n > \rho \log n)
    &\le \sum_{x \in B_n} P_p (\mathrm{diam}(\mathcal{C}_x) > \rho \log n, |\mathcal{C}_x|<\infty)\\
    &\le L n^d (\rho \log n)^d \exp\{-\xi(p)\rho\log n\}\\
    &= L \exp\left\{ d \log n + d \log (\rho \log n) - \rho\xi(p)\log n \right\}\\
    &= L \exp\left\{ \left( d - \rho\xi(p) \right) \log n + d \log \rho + d \log(\log n) \right\},
\end{align*}
which establishes the upper bound in the statement of the theorem.

Next, to prove the corresponding lower bound, we use that the events
\[
\{\mathrm{diam}(\mathcal{C}_x) \le \rho\log n\text{ or }\mathrm{diam}( \mathcal{C}_x) > 2 \rho\log n \}_{x\in A_n}
\]
are independent when
\begin{align*}
    A_n = (K_n \mathbb{Z})^d \cap B_n,
\end{align*}
and
\begin{align*}
    K_n = \lceil 2 (2\rho + 1) \log n \rceil . 
\end{align*}
By independence, we obtain
\begin{align*}
    P_p(R_n \leq \rho \log n )
    &\le  P_p \left( \bigcap_{x \in A_n}   \bigl\{ \mathrm{diam}(\mathcal{C}_x)\le \rho\log n\text{ or }\mathrm{diam}(\mathcal{C}_x)>2\rho\log n  \bigr\} \right)\\
    &=  \prod_{x \in A_n} \bigl(1-P_p(\mathrm{diam}(\mathcal{C}_x)\in (\rho\log n,2\rho\log n])\bigr)\\
    &=  \bigl(1-P_p(\mathrm{diam}(\mathcal{C}_0)\in (\rho\log n,2\rho\log n])\bigr)^{|A_n|}\\
    &\leq 1- \frac{1}2 |A_n|P_p(\mathrm{diam}(\mathcal{C}_0)\in (\rho\log n,2\rho\log n]),
\end{align*}
where we use the approximation $(1 - x)^m \le 1-xm/2$ for $m\in \mathbb N$ and  $x > 0$ as long as $xm$ is sufficiently small.
Moreover,  we use the fact 
\begin{align*}
    |A_n| \ge \left( \frac{n}{ 8 \rho\log n}\right)^d,
\end{align*}
so we arrive at a bound, for any $\delta>0$,
\begin{align*}
    P_p(R_n > \rho \log n )
    &\ge \frac 12 \left(\frac{n}{8 \rho\log n}\right)^d P_p(\mathrm{diam}(\mathcal{C}_0)\in (\rho\log n,2\rho\log n])\\
    &\ge \frac 12  \left(\frac{n}{8 \rho\log n}\right)^d \exp \bigl( -(\xi (p)+\delta)\rho\log n \bigr)\\
    &= \frac 12 \exp \left( \left( d - \rho(\xi(p)+\delta) \right) \log n - d \log (8 \rho) - d \log (\log n) \right),
\end{align*}
which completes the proof.
\end{proof}

\begin{proof}[\textnormal{\textbf{Proof of Theorem \ref{thm:ProbCon:Sn}}}]
First, we prove the following two-sided bound on the expectation.
For $x\in B_n$, let $D_x$ be the event that  $\{\rho \log n < \mathrm{diam}(\mathcal{C}_x) < \infty \}$. Then, 
the quantity $S_n:=S_n (\rho)$ can be described as the total number of such vertices in $B_n$: 
\[
S_n (\rho) = \sum_{x \in B_n} 1_{D_x}.
\]
By linearity, we have
\begin{align*}
E_p[S_n (\rho)] 
&= E_p\left[ \sum_{x \in B_n} 1_{D_x} \right] = \sum_{x \in B_n} E_p[1_{D_x}]= \sum_{x \in B_n} P_p (D_x).
\end{align*}
Since $P_p(D_x) = P_p(D_0)$ for any $x \in B_n$ by translational invariance, we obtain
\begin{align} \label{eq:Exv_sn}
E_p[S_n (\rho)] = |B_n| P_p(D_0).
\end{align}
By \eqref{def:xi'}, for any $\varepsilon > 0$ and sufficiently large $n$, $P_p (D_0)$ is sharply bounded as
\[
n^{- \xi(p)\rho - \varepsilon} \le P_p (D_0) \le n^{- \xi(p)\rho + \varepsilon}.
\]
Since $|B_n| = (2n+1)^d$, we have $n^d \le |B_n| \le (3n)^d$ for $n\ge 1$. Combining this with \eqref{eq:Exv_sn} yields the following bounds for the expectation:
\begin{align}\label{eq: expectation number}
    n^{d - \xi(p)\rho - \varepsilon} \le E_p[S_n (\rho)] \le 3^d n^{d - \xi(p)\rho + \varepsilon}.
\end{align}

    Next, we prove that $S_n (\rho)/E_p[S_n (\rho)]$ converges to $1$ in probability.
We aim to show that the right-hand side of the Chebyshev's inequality tends to $0$ as $n \to \infty$.
For any $\delta > 0$, the inequality is given by
\begin{align} \label{Chebyshev:Sn} 
P_p\left( \left|\frac{S_n (\rho)}{E_p[S_n(\rho)]}-1 \right|\ge \delta \right) \le \frac{\operatorname{Var}[S_n(\rho)]}{(\delta E_p[S_n(\rho)])^2}.
\end{align}
We write $I_x = 1_{D_x}$.
The variance is expanded as
\begin{align} \label{Var:var2cov} 
\operatorname{Var}[S_n] = \sum_{x \in B_n} \operatorname{Var} [I_x] + \sum_{x \in B_n} \sum_{y \in B_n : y \neq x} \operatorname{Cov} [I_x, I_y]. 
\end{align}
We proceed by estimating $\operatorname{Var} [I_x]$ and $\operatorname{Cov} [I_x, I_y]$ separately.

First, we estimate $\operatorname{Var} [I_x]$.
Since $I_x$ is an indicator function, we have that
\begin{align*} 
\operatorname{Var} [I_x] 
&= E_p[I_x^2] - E_p[I_x]^2\\
&= E_p[I_x] - E_p[I_x]^2\\
&= \bigl( 1 - P_p (D_x) \bigr)P_p (D_x)
\le P_p (D_x). 
\end{align*}
Summing over all $x \in B_n$, we obtain the following upper bound:
\begin{align} \label{Var:esn} 
\sum_{x \in B_n} \operatorname{Var} [I_x] 
&\le \sum_{x \in B_n} P_p (D_x) = E_p[S_n]. 
\end{align}

Next, letting $\eta := \max\{\varkappa(p), 2 \rho - \varkappa(p)\}$, we estimate $\operatorname{Cov}[I_x, I_y]$ by splitting $B_n$ based on the distance threshold $3\lceil (\varkappa(p) + \eta) \log n +2\rceil$ as follows:
\[
\sum_{y \in B_n : y \neq x} \operatorname{Cov}[I_x, I_y] = \sum_{
y \in x + B_{3\lceil (\varkappa(p) + \eta) \log n + 2\rceil},\,
y \neq x } \operatorname{Cov}[I_x, I_y] + \sum_{y \notin x + B_{3\lceil (\varkappa(p) + \eta) \log n + 2\rceil} } \operatorname{Cov}[I_x, I_y].
\]

First, we estimate for $y \in x + B_{3\lceil (\varkappa(p) + \eta) \log n + 2\rceil}$.
In this case, we apply the bound derived directly from the definition of covariance.
\begin{align*}
\operatorname{Cov}[I_x, I_y] 
&= E_p[I_x I_y] - E_p[I_x]E_p[I_y] \\
&= P_p(D_x \cap D_y) - P_p(D_x)P_p(D_y) \\
&\le P_p(D_x \cap D_y) \\
&\le P_p(D_x).
\end{align*}
Thus, the first sum is bounded as follows:
\begin{align}\label{estimate:shortrange}
\sum_{y \in x + B_{3\lceil (\varkappa(p) + \eta) \log n + 2\rceil},\,
y \neq x} \operatorname{Cov}[I_x, I_y]
&\le (7(\varkappa(p) + \eta) \log n)^d P_p(D_x),
\end{align}
where we use $|x + B_{3\lceil (\varkappa(p) + \eta) \log n + 2\rceil}| \le (7(\varkappa(p) + \eta) \log n)^d$ for sufficiently large $n$. 

Next, we estimate for $y \notin x + B_{3\lceil (\varkappa(p) + \eta) \log n + 2\rceil}$.
Recall from Theorem \ref{thm:LLN:Rn} and \eqref{ineq:Pp(diam)<=Ke^-xi} that for sufficiently large $n$, the probability that a cluster extends beyond the radius $\lceil (\varkappa(p) + \eta) \log n + 2\rceil$ is bounded as follows:
\begin{align} \label{ineq:diam_decay_explicit}
P_p(\mathrm{diam}(\mathcal{C}_x) > \lceil (\varkappa(p) + \eta) \log n + 2\rceil, |\mathcal{C}_x| < \infty) 
&\le L (\lceil (\varkappa(p) + \eta) \log n + 2\rceil)^d n^{-\xi(p)(\varkappa(p) + \eta)} \notag \\
&\le L' (\log n)^d n^{-d - \eta\xi(p)},
\end{align}
with some constants $L,L'>0$. We define the event $D'_x$:
\[
D'_x := D_x \cap \{ \mathcal{C}_x \subset x + B_{\lceil (\varkappa(p) + \eta) \log n + 2 \rceil} \}.
\]
The probability of $D_x \setminus D'_x$ satisfies \eqref{ineq:diam_decay_explicit}:
\begin{align} \label{bound:error_prob}
P_p(D_x \setminus D'_x) \le L' (\log n)^d n^{-d - \xi(p)\eta}.
\end{align}
We write $I'_x := 1_{D'_x}$ and decompose the indicator as $I_x = I'_x + 1_{D_x \setminus D'_x}$. Using the bilinearity of covariance, we expand $\operatorname{Cov}[I_x, I_y]$:
\begin{align*}
\operatorname{Cov}[I_x, I_y]
&= \operatorname{Cov}[I'_x + 1_{D_x \setminus D'_x}, I'_y + 1_{D_y \setminus D'_y}] \\
&= \operatorname{Cov}[I'_x, I'_y] + \operatorname{Cov}[I'_x, 1_{D_y \setminus D'_y}] + \operatorname{Cov}[1_{D_x \setminus D'_x}, I'_y] + \operatorname{Cov}[1_{D_x \setminus D'_x}, 1_{D_y \setminus D'_y}].
\end{align*}
Using the general property that for any indicator variables $X$ and $Y$, $\operatorname{Cov}[X,Y] \le E_p[X]$ holds, we obtain the bound:
\begin{align*}
\operatorname{Cov}[I_x, I_y]
&\le |\operatorname{Cov}[I'_x, I'_y]| + E_p[1_{D_y \setminus D'_y}] + E_p[1_{D_x \setminus D'_x}] + E_p[1_{D_x \setminus D'_x}] \\
&= |\operatorname{Cov}[I'_x, I'_y]| + 3 P_p(D_x \setminus D'_x),
\end{align*}
where we use the translational invariance $P_p(D_x \setminus D'_x) = P_p(D_y \setminus D'_y)$.
Since $y \notin x + B_{3\lceil (\varkappa(p) + \eta) \log n +2\rceil}$, the sets of edges on which $D'_x$ and $D'_y$ depend are disjoint. Hence $I'_x$ and $I'_y$ are independent and $\operatorname{Cov}[I'_x, I'_y]=0$.
Substituting this and \eqref{bound:error_prob} into the inequality, and summing over $y$, we obtain:
\begin{align} \label{estimate:longrange}
\sum_{y \notin x + B_{3\lceil (\varkappa(p) + \eta) \log n + 2\rceil}} \operatorname{Cov}[I_x, I_y]
&\le \sum_{y \notin x + B_{3\lceil (\varkappa(p) + \eta) \log n + 2\rceil}} 3 L' (\log n)^d n^{-d - \xi(p) \eta} \notag \\
&\le (3n)^d 3 L' (\log n)^d n^{-d - \xi(p)\eta} \notag \\ 
&=  3^{d+1} L' (\log n)^d n^{-\xi(p)\eta}.
\end{align}

By \eqref{estimate:shortrange} and \eqref{estimate:longrange}, we evaluate the total variance:
\begin{align*}
\operatorname{Var}[S_n] 
&\le E_p[S_n] + \sum_{x \in B_n} \left( (7 (\varkappa(p) + \eta) \log n)^d P_p(D_x) + 3^{d+1} L' (\log n)^d n^{-\xi(p)\eta} \right) \\
&\le E_p[S_n] + L''(\log n)^d E_p[S_n] + L''' (\log n)^d n^{d-\xi(p)\eta},
\end{align*}
for some constants $L'',L'''>0$.

Finally, we check the convergence condition:
\begin{align*}
\frac{\operatorname{Var}[S_n]}{(\delta E_p[S_n])^2} 
&\le \frac{E_p[S_n] + L'' (\log n)^d E_p[S_n] + L''' (\log n)^d n^{d-\xi(p)\eta}}{(\delta E_p[S_n])^2}\\
&\le \frac{1 + L'' (\log n)^d}{\delta^2 E_p[S_n]}
   + \frac{L''' (\log n)^d n^{d-\xi(p)\eta}}{\delta^2\, (E_p[S_n])^2}.
\end{align*}
By \eqref{eq: expectation number}, for any $\varepsilon>0$ and all sufficiently large $n$,
$E_p[S_n]\ge n^{d - \xi(p)\rho -\varepsilon}$, we obtain
\begin{align*}
\frac{1 + L'' (\log n)^d}{\delta^2 E_p[S_n]}
&\le \frac{1 + L'' (\log n)^d}{\delta^2}\,n^{-(d - \xi(p)\rho) +\varepsilon},
\end{align*} 
since $d - \xi(p)\rho > 0$ (as $\rho < \varkappa(p)$), taking $\varepsilon$ sufficiently small ensures the exponent is negative, so the first term converges to $0$ as $n \to \infty$.

For the second term, using $E_p[S_n]^2 \ge n^{2(d - \xi(p)\rho) - 2\varepsilon}$, we have
\begin{align*}
\frac{L''' (\log n)^d n^{d-\xi(p)\eta}}{\delta^2\, (E_p[S_n])^2}
&\le \frac{L'''}{\delta^2} (\log n)^d n^{d-\xi(p)\eta}\,n^{-2(d - \xi(p)\rho)+2\varepsilon}.
\end{align*}
The exponent of $n$ simplifies as follows: if we take $\varepsilon$ sufficiently small then, 
\begin{align*}
d-\xi(p)\eta - 2(d - \xi(p)\rho) + 2\varepsilon
&= -d + 2\xi(p)\rho - \xi(p)\eta + 2\varepsilon < 0,
\end{align*}
with the choice of $\eta = \max \{ \varkappa(p), 2 \rho - \varkappa(p) \}$. Hence,  we have
\[
(\log n)^d n^{-d + 2\xi(p)\rho - \xi(p)\eta + 2\varepsilon} \longrightarrow 0
\quad\text{as }n\to\infty,
\]
and hence the second term also converges to $0$.

Combining the above estimates, we conclude that
\[
\lim_{n\to\infty}
\frac{\operatorname{Var}[S_n]}{(\delta E_p[S_n])^2}
= 0,
\]
and therefore
\[
\frac{S_n(\rho)}{E_p[S_n(\rho)]} \xrightarrow{P_p} 1
\quad\text{as }n\to\infty.
\]
This completes the proof.
\end{proof}

\bigskip
\noindent \textbf{Acknowledgments.}\ 
The author would like to thank Prof. Shuta Nakajima (Keio University) for his continuous encouragement and insightful comments on this work.
The author is also grateful to Prof. Shigetoshi Yazaki (Meiji University) for providing valuable suggestions throughout this research.

\smallskip
\noindent \textbf{Kaito Kobayashi:}\ 

\smallskip \noindent
Department of Mathematics, Graduate School of Science and Technology, Meiji University,
1-1-1 Higashi-Mita, Tama-ku, Kawasaki, Kanagawa 214-8571, Japan.

\noindent E-mail: \texttt{ce246007@meiji.ac.jp}

\begin{thebibliography}{99}

\bibitem{BH1957}
S.~R.~Broadbent and J.~M.~Hammersley,
Percolation processes: I. Crystals and mazes,
\textit{Math. Proc. Cambridge Philos. Soc.}
\textbf{53} (1957), no.~3, 629--641.

\bibitem{SCA1982}
D.~Stauffer, A.~Coniglio and M.~Adam,
Gelation and critical phenomena,
in \textit{Polymer Networks} (K.~Dušek, ed.),
Advances in Polymer Science, vol.~44,
Springer, Berlin--Heidelberg, 1982.

\bibitem{SS}
M.~Sahimi,
\textit{Applications of Percolation Theory}, 2nd ed.,
Applied Mathematical Sciences, vol.~213,
Springer, Cham, 2023.

\bibitem{GGP2}
G.~Grimmett,
\textit{Percolation}, 2nd ed.,
Springer, Berlin, 1999.

\bibitem{KS1978}
H.~Kunz and B.~Souillard,
Essential singularity in percolation problems and asymptotic behavior of cluster size distribution,
\textit{J. Stat. Phys.} \textbf{19} (1978), 77--106.

\bibitem{ADS1980}
M.~Aizenman, F.~Delyon and B.~Souillard,
Lower bounds on the cluster size distribution,
\textit{J. Stat. Phys.} \textbf{23} (1980), 267--280.

\bibitem{JMH1957}
J.~M.~Hammersley,
Percolation processes: Lower bounds for the critical probability,
\textit{Ann. Math. Statist.}
\textbf{28} (1957), no.~3, 790--795.

\bibitem{K1987}
H.~Kesten,
Scaling relations for $2$D-percolation,
\textit{Commun. Math. Phys.}
\textbf{109} (1987), no.~1, 109--156.

\bibitem{HvH2007}
M.~Heydenreich and R.~van~der~Hofstad,  
Random graph asymptotics on high-dimensional tori,  
\textit{Commun. Math. Phys.}  
\textbf{270} (2007), no.~2, 335--358.   

\bibitem{HK2010}
B.~M.~Hambly and T.~Kumagai,  
Diffusion on the scaling limit of the critical percolation cluster in the diamond hierarchical lattice,  
\textit{Commun. Math. Phys.}  
\textbf{295} (2010), no.~1, 29--69.  

\bibitem{KB2016}
R.~Kenna and B.~Berche,  
Universal finite-size scaling for percolation theory in high dimensions,  
\textit{J. Phys. A: Math. Theor.}  
\textbf{50} (2017), no.~23, 235001.  

\bibitem{AKN1987}
M.~Aizenman, H.~Kesten and C.~M.~Newman,
Uniqueness of the infinite cluster and continuity of connectivity functions for short and long range percolation,
\textit{Commun. Math. Phys.}
\textbf{111} (1987), 505--531.

\bibitem{BK1989}
R.~M.~Burton and M.~Keane,
Density and uniqueness in percolation,
\textit{Commun. Math. Phys.}
\textbf{121} (1989), 501--505.

\bibitem{KESTEN1980}
H.~Kesten,
The critical probability of bond percolation on the square lattice equals \(1/2\),
\textit{Commun. Math. Phys.}
\textbf{74} (1980), 41--59.

\bibitem{AGKSZ2014}
N.~Araújo, P.~Grassberger, B.~Kahng, K.~J.~Schrenk and R.~M.~Ziff,
Recent advances and open challenges in percolation,
\textit{Eur. Phys. J. Special Topics}
\textbf{223} (2014), 2307--2321.

\bibitem{DC2018}
H.~Duminil-Copin,
Sixty years of percolation,
in \textit{Proc. ICM}, Vol.~IV,
World Scientific, Hackensack, NJ, 2018, 2829--2856.

\bibitem{G2023}
G.~Grimmett,
Selected problems in probability theory,
in \textit{Lecture Notes in Math.} 2313,
Springer, Cham, 2023, 603--614.

\bibitem{M2017}
R.~Morris,
Bootstrap percolation, and other automata,
\textit{European J. Combin.}
\textbf{66} (2017), 250--263.

\bibitem{BS1996}
I.~Benjamini and O.~Schramm,
Percolation beyond $\mathbb{Z}^d$, many questions and a few answers,
\textit{Electron. Commun. Probab.}
\textbf{1} (1996), 71--82.

\bibitem{M1986}
M.~V.~Menshikov,
Coincidence of critical points in percolation problems,
\textit{Soviet Mathematics Doklady}
\textbf{33} (1986), 856--859.

\bibitem{AB1987}
M.~Aizenman and D.~J.~Barsky,
Sharpness of the phase transition in percolation models,
\textit{Commun. Math. Phys.}
\textbf{108} (1987), no.~3, 489--526.

\bibitem{CCGKS1989}
J.~T.~Chayes, L.~Chayes, G.~R.~Grimmett, H.~Kesten and R.~H.~Schonmann,
The correlation length for the high-density phase of Bernoulli percolation,
\textit{Ann. Probab.}
\textbf{17} (1989), no.~3, 1277--1302.

\bibitem{DN2022}
B.~Dembin and S.~Nakajima,
On the upper tail large deviation rate function for chemical distance in supercritical percolation,
arXiv:2211.02605.

\bibitem{HR2006} 
R.~van~der~Hofstad and F.~Redig, Maximal clusters in non-critical percolation and related models, \textit{J. Stat. Phys.} 
\textbf{122} (2006), 671--703.
\end{thebibliography}
\end{document}